\documentclass{amsart}

\title[An American  in Mathematician in Moscow]
{An American Mathematician in Moscow,
or How I Destroyed the Soviet Union}
\author{Melvyn B. Nathanson}
\address{Department of Mathematics,
Lehman College (CUNY),
Bronx, New York 10468, and 
CUNY Graduate Center, New York, NY 10016}
\email{melvyn.nathanson@lehman.cuny.edu}

\begin{document}

\maketitle

\begin{center}
Dedicated to I. M. Gel'fand on the 100th anniversary of his birth.
\end{center}
\vspace{1cm}

 The great Soviet mathematician Israel Moiseevich Gel'fand 
was born on September 2, 1913 in Okny (later Krasni Okny, or Red Okny) 
near Odessa in the Ukraine, 
and died on October 5, 2009 in New Brunswick, New Jersey.  
The Russian Revolutions of 1917 led to the approval of the Treaty of Creation 
of the USSR on December 29, 1922.  
The Soviet Union ceased to exist on December 26, 1991.  
Gel'fand was born before the Soviet Union, and outlived it.  

I was indirectly introduced to Gel'fand in 1970.  
I was a visiting research student in the Department of Pure Mathematics 
and Mathematical Statistics  of the University of Cambridge 
during the Lent and Easter terms.   
One of my friends  was B\' ela Bollob\' as, a Hungarian who had received 
his Ph.D. at Oxford and had decided to remain in England and not return to Hungary.  
In the terminology of the Cold War, B\' ela had defected to the capitalists.  
Before being allowed to study in the West, the Hungarians had required 
him to  study in the East, that is, in the USSR, and B\' ela  had spent a year 
at Moscow State University, where he worked with  Gel'fand. 
Professor Gel'fand had impressed him deeply as a   man  and as a mathematician, 
and B\' ela often told me how extraordinary  he was.  

In the summer of 1970, at the end of my study in Cambridge, 
I made a short trip to the Soviet Union.  
Another American, a biologist from MIT,  
had just completed a post-doctoral year 
at Cambridge and had posted a note on a bulletin board that he was looking 
for someone to accompany him on a driving tour through the USSR.  
The plan was to enter the Soviet Union by train from Helsinki to Leningrad, 
rent a car, and drive south through Moscow to the Caucasus.  
The biologist had been invited to lecture in Moscow as a guest 
of the Soviet Academy of Sciences.  His official host was a distinguished 
Soviet biochemist, David Gold'farb.  I also met Gold'farb, and told him that 
I would like to visit Gel'fand.  Gold'farb contacted Gel'fand, who was too busy 
or (more likely) too prudent to rendezvous with an unknown American.  
When Gold'farb told Gel'fand that I was interested in number theory, 
Gel'fand gave him to give to me a copy of his book 
\emph{Representation Theory and Automorphic Functions}, 
written with M. I. Graev and I. I. Pyatetskii-Shapiro, 
the sixth volume of the series of monographs \emph{Generalized Functions}.  

Gold'farb had lost a leg fighting in World War II.  
He had wanted to be a historian, but history is a dangerous profession 
in totalitarian regimes.  
In Stalin's Russia,  history was particularly life-threatening, 
so Gold'farb went to medical school 
and did research in molecular genetics.  
``A stomach is always a stomach," he told me.

My biologist traveling companion was very left politically.  
Like most academics, I was against the war in Vietnam, 
but he was so far to my left that by comparison 
I seemed to be on the right.  
That immediately endeared to me to Soviet scientists, 
many of whom acted in public as if they were  loyal  followers 
of the Communist Party line, but inwardly were strongly anti-totalitarian 
and, indeed,  unlike most  American scientists, supported American intervention 
in Vietnam.   They believed in killing Commies.  
This was one of the first lessons in irony that I learned in Moscow.  

Many senior Communist Party officials pulled strings to get their children into scientific careers.  They had done what they needed to do to survive, 
they understood the dangers of politics, and they wanted safer lives 
for their own children.  
Many graduate students and researchers in the Soviet Academy of Sciences and in Soviet universities were relatives of high-ranking Party functionaries.  
Of course, being a scientist with ties to the Communist Party brought its own  benefits.  
Many distinguished non-Party Soviet researchers collaborated with colleagues 
who were much weaker scientists, but politically powerful,  
in order to shelter themselves and their students.  
``Collaboration'' often meant nothing more than adding the names of mediocre 
Party members  to their papers as co-authors.  
Of course, this kind of scientific politics also occurs in non-totalitarian regimes.   

It was because of B\' ela Bollob\' as that I got the idea of trying 
to spend a year studying mathematics with Gel'fand. 
This would not be easy to arrange.  
During the Cold War it was almost impossible for an American to study 
or do research in Moscow.  To study in Cambridge or Paris was trivial.  
Just get on a plane and fly to England or France.  
But the only way an ordinary American could enter the Soviet Union 
was on a brief and expensive tourist visit.  There were, however, two formal academic 
exchanges.  One was between the National Academy of Sciences in Washington 
and the Soviet Academy of Science, but this usually provided only short-term visits 
for senior scientists, not young scholars.  

The other program was part of a broad cultural affairs treaty 
between the US and the USSR.  We would send the New York Philharmonic 
to Moscow, and they would send the Bolshoi Ballet to New York.  
One small part of this treaty was a university-level exchange 
for graduate students and post-docs.  
Each year the Americans selected 40 American scholars, and
the Soviets selected 40 Soviet scholars.  
Each country's choices had to be approved by the  other.  
The US program was  administered by IREX, the International Research 
and Exchanges Board, an organization based in New York 
and associated with the American Council of Learned Societies.    

Typically, the Russians sent 40 engineers and computer scientists to MIT, 
and the Americans sent 40 specialists in literature and history to Moscow.  
IREX had sent very few scientists to the USSR.   
The logic on the American side was reasonable:  
Because there were so few opportunities to do research in the USSR, 
an American exchange student should have a research project 
that could not be carried out 
anywhere else in the world.  Science is everywhere, so it would be hard to argue 
that a scientific research problem could only be solved in the USSR.  
On the other hand, if your scholarly work were in Russian or Soviet  
literature or history, 
and if the archives you needed were in the Soviet Union, 
then you clearly had a research 
proposal that would justify a trip to Moscow.  

In 1972-73 I became  the first American mathematician to participate 
in the IREX exchange program.  
In my application to IREX, I wrote that the Soviet Union had many 
of the greatest mathematicians in the world, that they were concentrated 
in Moscow and not allowed to travel outside the country, and that it would be 
extremely valuable to be able to meet and work with them.  
The argument had merit, and was accepted by IREX.
My research proposal was to work with Gel'fand.  
He had to agree to supervise me, and he did.  
I got my visa and went.

The Americans on the IREX exchange were supposed to rendezvous  
in Paris in August and fly together to Moscow.  
I got to Paris a few days early with a suitcase full of math books 
and an Olivetti portable typewriter, stayed in a cheap hotel on Rue des \' Ecoles, 
and worked desperately hard to finish what would become my first joint paper 
with Paul Erd\H os.  
I did not have time to mail the manuscript from France, but 
IREX exchange students had certain privileges at the American Embassy in Moscow, 
and one of the  most valuable was the use of  the diplomatic pouch 
for sending letters out of the USSR.   
The Soviet postal system was, to put it politely, ``unreliable'' for manuscripts 
being sent to the United States, but I was able to submit my paper 
to the  \emph{Proceedings of the AMS}  
in a mailsack hand-carried by a US Marine to Vienna.  

My first meeting with Gel'fand was in the lobby of Moscow State University.  
I remember two things that he told me.  
The first was his famous mantra:  ``There is only one mathematics.''
Then, after reciting a short list of the best young mathematicians 
in the Soviet Union, he said,  
``They know much more mathematics than I, but my intuition is better."
Gel'fand suggested that I attend courses by Pyatetskii-Shapiro and Manin, but the most important part of my mathematical education in Moscow was participation 
in Gel'fand's famous Monday night seminar.  
I don't recall the official starting time of the seminar.  
People would show up early and talk mathematics in the hall, the seminar would eventually begin, and there would be a series of speakers, lasting long after the seminar was supposed to end, until finally we were evicted by a cleaning lady who had to do her job.  
 
It was common in the seminar for Gel'fand to interrupt a talk and ask someone 
in the audience to explain what was going on.  The first time I went to the seminar, 
in the middle of a lecture, he asked, "Melvyn, do you understand?"

"No," I replied.  

``Why not?"  

``Because they're speaking Russian."

He then assigned Dima Fuks the task of sitting next to me and translating the lecture from Russian into English.  In a short time my Russian improved enough that I could understand the talks, and the language excuse was lost.  

Gel'fand would decide that someone needed to learn something 
and present it in the seminar.  
For example, he asked Arnol'd to give a series of talks on $p$-adic numbers.  
Arnol'd seemed to find this difficult.  A master of the real and complex domains, 
he had trouble understanding non-Archimedian absolute values.  
Of course, this is something every young number theorist knows.  
To see a great mathematician like Arnol'd struggling with $p$-adic analysis  
teaches that you are not an idiot if you don't understand some piece 
of mathematics  that ``everyone'' finds trivial.  

Walking is a  Russian tradition.  
The winters are cold, but there is little wind and the effective temperature is 
certainly bearable.  After the seminars, a group of people would often leave 
with Gel'fand, and walk and talk late into the night, 
while writing mathematics in the snowdrifts along the sidewalks.  
Outside, you could talk more freely than in rooms where the walls had ears.  
In the course of the year, many mathematicians would ask me to go for a walk, 
and in the privacy of the streets would ask, 
``What is it like in America?''
``How much anti-Semitism is in America?''
``How hard is it to get a job in an American university?''  
In a few years, as soon as emigration became possible, they all emigrated.

Gel'fand immigrated to the United States in 1989.  He was a visiting professor 
at Harvard and MIT, and then distinguished professor at Rutgers University.  
The biochemist David Gold'farb also left Russia for New York.  

American students on the IREX exchange lived in the dormitory of Moscow State 
University.  We were told that Americans were always assigned the same rooms, 
not on the same floor, but on different floors, one room directly above the other, 
because it was easier to bug them by dropping wires vertically through the building.  
A standard joke:  ``If you need something in your room fixed, speak into the lightbulb."  
I had many friends who were active in the university Komosomol, 
the youth division of the Communist Party.  
They were, presumably, assigned to befriend Americans.  
One of them told me, ``They can't identify all the voices on the tapes from your room."  

The Komsomol mirrored life outside the university, where Communist Party leaders 
had perks and privileges not available to the hoi-poloi.  
For example, there were private parties in the university for the Komsomol elite only.   
I attended a party with music provided by a rock band brought in from Estonia.  
When ``ordinary students'' tried to crash the Party, the Komsomol called the police.

I flew back to Philadelphia during the Christmas break to visit my mother, who was sick.  
When I returned, Gel'fand asked, ``How did it feel to be in the Soviet Union, 
then back in the US, then back in Moscow?"  I replied, "OK, but for the first week 
at home I was afraid to use the phone."

Gel'fand urged me to read Russian literature, especially Pushkin.  
He gave me a recording of Pushkin's poem, \emph{Mozart and Salieri}, 
and novels by Ilf and Petrov.  He also gave me various mathematics books, 
including his book with Minlos and Shapiro, 
\emph{Representations of the Rotation Group and Lorentz Group}, 
and the English edition of Weil's \emph{Basic Number Theory}.  Gel'fand had an enormous capacity to create friendships.  
Andr\' e Weil visited him in Moscow, and they became close.  
After I returned to the US, Weil invited me to spend a year as his Assistant 
at the Institute for Advanced Study. 
I do not know, but always assumed, that Gel'fand had recommended me, 
and his friend Weil obliged.  

You learned in Moscow to keep your Soviet friends in disjoint circles.
Knowing an American was dangerous, informants were ubiquitous, 
and some of your acquaintances  were undoubtedly reporting on you 
to the ``competent organs," 
which wanted to know  everyone with whom an American was in contact.    
There was no reason for me to be paranoid, only careful.  
I always felt   completely safe because I held an American passport.  
The Russians would not want to create an international incident.  
I might be arrested and threatened, but if I kept cool I would only be deported, 
which was no big deal.  But Soviet citizens could really be endangered, 
expelled from universities, fired from jobs, their lives seriously impacted.  
Even though the United States was intensely waging the Vietnam war 
and we were bombing Hanoi, 
Kissinger's policy of a multi-track foreign policy with the stick in Southeast Asia
 and the carrot to the USSR, namely, the allure of American trade concessions 
 and exports  to Russia, convinced me that an American in Moscow 
 on an official academic exchange 
 program whose only crime was  ``acting like an American,'' not espionage, 
 was perfectly safe.

In Moscow State University, as in all  Soviet universities, 
there was the ``First Department'' (in Russian, the ``Pervii Otdel''), 
which was the KGB office within the university.  
After  Gel'fand had emigrated and was a professor at Rutgers, 
he recounted the following story.
``I could not tell you this when you were in Moscow," he said, ``but during your stay here I was visited by someone from the Pervii Odtel.   
The KGB officer told me, `You have an American student, Nathanson.   
I know Americans are independent, but Nathanson is too independent 
even for an American, and we have to expel him from the country.' "

Gel'fand, who possessed great political savvy, replied, 
``Of course, you should expel him if you have to, but  I know that Nathanson 
has many important friends in America, and, if you expel him, 
there will be an international furor.  It might be better to let him finish the year, 
and then not let him return."  
That's what happened.  

Gel'fand thought it would be good for my mathematical education to spend another year 
with him in Moscow, but it would clearly be impossible for me to return 
to Moscow State University.  
The other US-USSR scientific exchange program was with the Soviet Academy of Sciences.  
Gel'fand told me to apply to the Academy exchange, but not to request placement 
in the Stekhlov Institute, which was the notoriously anti-Semitic 
mathematics  institute in the Soviet Academy.  
Instead, Gel'fand recommended that I apply to the 
Institute for Problems in the Transmission of Information, 
where several first-rate Soviet Jewish mathematicians found  safe haven.  
In 1977 I applied and was accepted, but at the last minute 
the Soviet Foreign Ministry refused to issue me a visa and I could not go.  
This action did, in fact, become an international incident, 
with coverage in \emph{The New York Times} and news media around the world, 
as a Soviet violation 
of the human rights for scientists provisions of the Helsinki Accords.

Under Communism, whether in the Soviet Union or in Eastern Bloc countries, 
there was a strange and tense separation of one's inner life and outer life, 
between what one had to say and do in front of strangers, 
and how one thought and acted with friends you really trusted.  
It was, as Russians liked to say, ``sloznii,'' that is, ``complicated.''   
Academic jobs in in Moscow typically  went to those 
who were well connected and acceptable to the Communist Party.  
There was no great monetary reward for studying mathematics, or, more precisely, 
for living mathematics.   It was done for free, for love, for intellectual 
and emotional enrichment, and not, 
as often  in the West, for professional advancement.  
With the collapse of the USSR, Russian mathematics lost some of its purity 
and became more, in the American and European sense, ``professional.''  
One other hand, now you can buy meat in Irkutsk, 
so we in the West, who never experienced Soviet-scale deprivation, 
should not be disparaging about this.  
To an American in Moscow  during the Cold War, the intellectual quality of  life in mathematical circles was awesome.  

Kazhdan once said that, when he got to know me,  
he had never met anyone with my attitude, a kind of unfrightened, 
relaxed approach to life.  In Russia everyone was constantly on guard, 
alert to danger, afraid of saying something that could unintentionally, 
or malignantly intentionally, be misinterpreted, reported to the ``competent authorities,'' 
and cause expulsion from school, exile, imprisonment, or death.  
One had to be always vigilant.  One had to make decisions:  Who do I trust?  
How much can I trust this person?  How open can I be?  
Is this guy reporting on me to the KGB?  
Will this person lie about me for no obvious reason?  
Americans don't understand this pressure.   We have grown up unpersecuted 
and without fear of persecution.  
Our country was rich, even if we were not, and there was a sense of fairness.  
Even an older generation, in the bad but brief period of Communist witch hunts and McCarthyism, never had to fear what Soviet citizens feared.  

It would be hard to overestimate the brittleness of the former USSR.  
An American could endanger its  political system 
by going to Moscow and being American.  
By not being terrified.  Not being cowed.  
Not being blackmailed by the threat of exclusion from libraries and archives.  
Soviets could sense the huge psychological difference between living in a free country 
and living in a totalitarian one.   Soviet authorities were correct to want to keep Americans 
away from ordinary Russians.  
We were a threat to the state.  
We were dangerous.

\end{document}